\documentclass[runningheads]{llncs}
\usepackage[T1]{fontenc}
\usepackage   	{algorithm,algorithmic}
\usepackage   	{amsmath}
\usepackage	  	{float}
\usepackage		{amssymb}
\usepackage		{hyperref}
\usepackage		{graphicx}
\usepackage		{color}

\algsetup{indent=2em}

\begin{document}
\title{A Matheuristic Multi-Start Algorithm for a Novel Static Repositioning Problem in Public Bike-Sharing Systems.}
\titlerunning{A Matheuristic Multi-Start Algorithm for SBRP}
\author{Julio Mario Daza-Escorcia\orcidID{0000-0003-4557-5150} \\
David Álvarez-Martínez\orcidID{0000-0001-8411-1936}} 
\authorrunning{Daza-Escorcia and Álvarez-Martínez}
\institute{Universidad de los Andes 
	\\ Cra. 1 No 18a - 12 Edif. Mario Laserna Pinzón (Bogotá), Colombia
	\url{https://pylopre.uniandes.edu.co/} \\
	Grupo PyLO Producción y Logística \\
	\email{\{jm.dazae,d.alvarezm\}@uniandes.edu.co}}
\maketitle

\begin{abstract}
This paper investigates a specific instance of the static repositioning problem within station-based bike-sharing systems. Our study incorporates operational and damaged bikes, a heterogeneous fleet, and multiple visits between stations and the depot. The objective is to minimize the weighted sum of the deviation from the target number of bikes for each station, the number of damaged bikes not removed, and the total time used by vehicles. To solve this problem, we propose a matheuristic approach based on a randomized multi-start algorithm integrated with an integer programming model for optimizing the number of operatives and damaged bikes that will be moved between stations and/or the depot (loading instructions). The algorithm's effectiveness was assessed using instances derived from real-world data, yielding encouraging results. Furthermore, we adapted our algorithm to a simpler problem studied in the literature, achieving competitive outcomes compared to other existing methods. The experimental results in both scenarios demonstrate that this algorithm can generate high-quality solutions within a short computational time.


\end{abstract}
%
%
%

\section {Introduction}

Bike-sharing systems (BSS) originated over 50 years ago in Northern Europe, serving as a mobility facilitator and a supplement to traditional public transportation. Notable advantages include their minimal environmental impact and cost-effectiveness compared to other transportation modes. Despite these benefits, specific challenges and instances of dissatisfaction need attention. For instance, issues arise when a station needs more free anchors for bike parking, has insufficient available bikes for users, or when the available bikes are damaged, hindering the system's smooth operation.

These scenarios give rise to what is known as the \emph{Bike-sharing Repositioning Problem} (BRP), involving the redistribution or rebalancing of bikes within the system through a fleet of vehicles. This process entails the movement of operative bikes and/or damaged bikes between stations and the depot to align their inventory with a desired or target level.

Two primary classifications exist within public bike-sharing systems. The first, referred to as the \emph{Station-Based Bike Sharing} (SBBS) system, allows users to rent bikes at designated stations and return them either to the same station or any specified station upon completion of use \cite{Alvarez-Valdes2016}. The second classification is the \emph{Free-Floating Bike Sharing} (FFBS) system, where bikes can be picked up and dropped off at locations chosen by the users \cite{Pal2017}.

The literature identifies two types of repositioning strategies: dynamic and static. In the \emph{dynamic approach}, bikes are repositioned while the system is actively used, such as during the day. Instead, the \emph{static approach} involves repositioning during minimal from the bike's utilization on the system, typically at night.

Our research focuses on a static approach to the bike-sharing repositioning problem, particularly within the station-based model. We refer to this as the \emph{Station-Based Static Bike-sharing Repositioning Problem}. This study provides significant contributions to the typical \emph{Static Bike-sharing Repositioning Problem} (SBRP) in the following aspects:

\begin{enumerate}
	\item Introducing a novel instance of the static repositioning problem within station-based bike-sharing systems, incorporating additional factors to enhance realism. These include considerations for operational and damaged bikes, a heterogeneous fleet, and multiple visits between stations and the depot.
	
	\item Presenting an integer programming model to optimize the number of operative and damaged bikes to be moved between stations and/or the depot (loading instructions).
	
	\item Proposing an effective matheuristic based on a randomized multi-start algorithm integrated with an integer programming model to handle large instances.
\end{enumerate}

The remainder of the paper is organized as follows. Section \ref{sec:Literature_Review} provides the literature review. The description of the research problem proposed is presented in Section \ref{sec:Problem_Description}. Section \ref{sec:Solution_Methodology} describes the proposed solution methodology to solve the research problem. Section \ref{sec:Randomized_Multi-Start_Algorithm} presents a Randomized Multi-Start Algorithm, and Section \ref{sec:Optimal_Loading_Instructions} describes an exact method to optimize loading instructions of a previously constructed set of routes. Section \ref{sec:Computational_Experiments} introduces the set of instances used in this study and our computational experiments. Finally, conclusions and remarks are presented in Section \ref{sec:Concluding_Remarks}.
\section {Literature Review}				\label{sec:Literature_Review}
Much research has been published recently exploring different aspects of the BRP. This paper primarily examines the characteristics of SBRP, including its objectives and optimization methods. The characteristics of SBRP include the fleets, multiple or single visits to the stations or depot, and types of bikes for repositioning.

The literature encompasses works that explore either homogeneous or heterogeneous fleets of vehicles. Most of these studies employ a heterogeneous fleet of vehicles, which is also the case in our problem. Including a heterogeneous fleet increases the complexity of the problem but results in a more realistic scenario. The studies for the static case that employing a heterogeneous fleet include those by Alvarez-Valdes et al.~\cite{Alvarez-Valdes2016}, Casazza~\cite{Casazza2016}, Di Gaspero et al.~\cite{DiGaspero2013,DiGaspero2013a,DiGaspero2016}, Espegren et al.~\cite{Espegren2016}, Forma et al.~\cite{Forma2015}, Kinable~\cite{Kinable2016}, Papazek et at.~\cite{Papazek2013,Papazek2014}, Raidl et al.~\cite{Raidl2013}, Rainer-Harbach et al~\cite{Rainer-Harbach2013,Rainer-Harbach2015}, Raviv et al.~\cite{Raviv2013}, Schuijbroek et al.~\cite{Schuijbroek2017}, Du et al.~\cite{Du2020} and Wang and Szeto~\cite{Wang2018}.

Multiple visits to the stations can occur in two distinct manners: (1) when a station is open for visits by multiple vehicles but restricts repeated visits by the exact vehicle (similar to the \emph{Split Delivery Vehicle Routing Problem}, SDVRP), and (2) when stations permit multiple visits by the same or different vehicles. In this study, we specifically focus on this latter scenario because it reflects a common occurrence in real-world settings within BRPs. This feature is evident in the works of Alvarez-Valdes et al.~\cite{Alvarez-Valdes2016}, Casazza~\cite{Casazza2016}, Di Gaspero et al.~\cite{DiGaspero2013,DiGaspero2013a,DiGaspero2016}, Espegren et al.~\cite{Espegren2016}, Forma et al.~\cite{Forma2015}, Kloimullner and Raidl~\cite{Kloimullner2017}, Papazek et at.~\cite{Papazek2013,Papazek2014}, Raidl et al.~\cite{Raidl2013}, Rainer-Harbach et al~\cite{Rainer-Harbach2013,Rainer-Harbach2015}, Raviv et al.~\cite{Raviv2013}, Schuijbroek et al.~\cite{Schuijbroek2017}, Du et al.~\cite{Du2020} and Wang and Szeto~\cite{Wang2018}.

Regarding the types of bikes, most studies focus on a single type of bike (operative or usable bikes). Alvarez-Valdes et al.~\cite{Alvarez-Valdes2016} and Wang and Szeto~\cite{Wang2018} acknowledge the potential presence of damaged (unusable) bikes, suggesting their removal from the stations. In contrast, Li et al.~\cite{Li2016} introduces diverse types of bikes but does not address the possibility of having damaged ones. Lastly, Du et al.~\cite{Du2020} stand out as they consider the collection of damaged bicycles, but specifically in the context of \emph{Free-Floating Bike Sharing} (FFBS) systems.

Beyond the diverse characteristics, the objectives explored in the literature of BRPs exhibit considerable variation. Typically, these studies seek to optimize one or multiple performance metrics, such as travel cost, time or distance, the count of loading and unloading operations at stations, the absolute deviation from the target number of bikes at stations, and so forth. In this paper, we propose a multi-objective problem aimed at minimizing the weighted sum of three terms: the deviation from the target number of bikes for each station, the number of damaged bikes left unremoved, and the total time used by vehicles.

Regarding the solution methods for the SBRP, while the majority are heuristics, we encountered some tailored algorithms. Table \ref{tab:Solution_Methods} summarizes the heuristic solution methods for SBRP. In our study, we propose a metaheuristic approach, which combines mathematical optimization techniques with heuristic methods to address complex optimization problems, such as the SBRP proposed here.
\begin{table}[!htbp]
	\caption{Heuristic solution methods for SBRP.} \label{tab:Solution_Methods}
	\begin{center}
		\scriptsize	
		\begin{tabular}{|l|l|}
			\hline 
			\qquad \qquad \qquad Solution method &  \qquad \qquad Article \\
			\hline
			3-step heuristic & \cite{Forma2015} \\					
			9.5-approximation Algorithm & \cite{Benchimol2011}  \\					
			Ant Colony Optimization & \cite{DiGaspero2013}  \\					
			Chemical Reaction Optimization & \cite{Szeto2016}  \\					
			Cluster-first Route-second Heuristic & \cite{Schuijbroek2017} \\					
			Destroy and Repair Algorithm & \cite{DellAmico2016} \\					
			Genetic Algorithm & \cite{Kadri2016}, \cite{Lahoorpoor2019}, \cite{Du2020} \\					
			Greedy Randomized Adaptive Search Procedure &  \cite{Papazek2013},  \cite{Papazek2014},  \cite{Rainer-Harbach2015}  \\					
			Heuristic based on Minimum Cost Flow Problem & \cite{Alvarez-Valdes2016} \\					
			Hybrid Genetic Search & \cite{Li2016} \\					
			Iterated Local Search &  \cite{Cruz2017}, \cite{Tang2020} \\					
			Iterated Tabu Search &  \cite{Ho2014}  \\					
			Large Neighborhood Search &  \cite{DiGaspero2013a}, \cite{DiGaspero2016},  \cite{Ho2017}, \cite{Pal2017}  \\					
			Particle Swarm Optimization &  \cite{Zhang2018} \\					
			Path Relinking &  \cite{Papazek2014}  \\					
			Tabu Search &  \cite{Chemla2013a}  \\					
			Variable Neighborhood Descent & \cite{Cruz2017}, \cite{Papazek2013}, \cite{Papazek2014}, \cite{Pal2017} , \cite{Rainer-Harbach2013}, \cite{Rainer-Harbach2015}, \cite{Zhang2018}  \\					
			Variable Neighborhood Search &  \cite{Raidl2013},  \cite{Rainer-Harbach2013}, \cite{Rainer-Harbach2015}, \cite{Lv2020} \\
			\hline
		\end{tabular}
	\end{center}
\end{table}
As apparent from the preceding literature review, several studies have been undertaken on the SBRP. Nonetheless, based on our understanding, more research is still needed concerning static repositioning bike-sharing problems, particularly regarding the collection of damaged bikes. This gap is especially notable when considering the inclusion of a heterogeneous fleet and multiple visits within this problem domain.

This paper presents a novel instance of the repositioning bike-sharing problem with a static approach and within the station-based model. This problem incorporates operational and damaged bikes, a heterogeneous fleet, and multiple visits between stations and the depot. The objective is to minimize the weighted sum of the deviation from the target number of bikes for each station, the number of damaged bikes left unremoved, and the overall time used by vehicles. To solve this problem, we propose a matheuristic approach combining a randomized multi-start algorithm with an integer programming model to optimize the number of bikes relocated during each visit.

\section {Problem Description}  			\label{sec:Problem_Description}
The \emph{Station-Based Static Bike-sharing Repositioning Problem} (SSBRP) considering in this paper is defined on a complete directed graph $G_o=(V_o, A_o)$, where $V_o$ are the nodes set represented to the set of stations $V$ and the depot $O$ (i.e, $V_o=V\cup O$), and $A_o$ are the arcs representing the shortest paths (one for each pair of nodes). Each arc $(u,v)\in A_o$ with $u,v\in V_o$ has a cost $t_{(u,v)}$ representing the travel time between $u$ and $v$.

The depot $O\in V_o$ has a sufficient capacity $c_o$, and an initial inventory of operative (or usable) bikes $p_o\geq 0$. Each station $v_\in V$ has a capacity (or parking docks) $c_v>0$, an initial inventory of operative bikes $p_v\geq 0$ and damaged bikes $a_v\geq 0$, a target (desired) ending inventory $q_v\geq 0$, and a weight (or visit priority) $w_v$.

We define the initial unbalance of operative bikes $d_v=p_v-q_v$. According to its imbalance, a station can be: balanced $V_{bal}=\{v\in V\mid d_v=0\}$, in surplus $V_{pic}=\{v\in V\mid d_v > 0\}$, in deficit $V_{del}=\{v\in V\mid d_v < 0\}$. In either case, a station may have damaged bikes $V_{ave}=\{v\in V\mid a_v\geq 0\}$.

A set of trucks $L$ can pick up or deliver bikes at each station or depot (i.e., $L=\{1,\dots,\lvert L\rvert\}$). Each vehicle $l\in L$ has a heterogeneous capacity of $k^{l}$, a route time no longer than $T$, and the depot as the starting and ending point of the tour.

A possible solution to the SSBRP consists of a set of $R$ routes and $Y$ movements. Each route $r^{l}\in R$ is assigned to a vehicle $l$, and has an ordered list of tours $r^{l}=\{v_1^{l},\ldots,v_{n^{l}}^l\}$. The route must start and end at the depot, and not exceed the capacity $k^{l}$ and the maximum route time $T$.

A loading instructions $y^{l}$ consists of a number of operative $B_i^{l}$ and damaged $A_i^{l}$ bikes that will be moved on visit $i$ by vehicle $l$, such that if $B_i^{l}>0$ operative bikes are loaded on the vehicle, on the other hand if $B_i^{l}<0$ these bikes are unloaded, if $B_i^{l}=0$ no movements of operative bikes are made. Similarly for $A_i^{l}$. The load instructions are formally defined as $y^{l}=\{(B_1^{l}, A_1^{l}),\ldots,(B_{n^{l}}^{l}, A_{n^{l}}^{l})\}$).

The proposed SSBRP has as its objective function to design the routes of the vehicles and the number of both operative and damaged bikes to be moved at each station in such a way that a weighted sum of three terms is minimized, where $\hat{d}_v$ and $\hat{a}_v$ represent the operative and damaged bikes at the end of the repositioning operation. $\gamma_d, \gamma_a, \gamma_t$ are the weights of the respective terms.
\begin{equation}
	\frac{\sum_{v \in V}w_v|q_v-\hat{p}_v|}{\sum_{v \in V}w_v|q_v-p_v|+a_v} \gamma_d +
	\frac{\sum_{v \in V}w_v\hat{a}_v}{\sum_{v \in V}w_v|q_v-p_v|+a_v} \gamma_a +
	\frac{\sum_{l \in L}t^l}{T*|L|}  \gamma_t 
	\label{FO-SBRP}
\end{equation}

The first and second terms represent the imbalance and the number of damaged bikes at the end of the repositioning for all the stations. The third term relates to the total fleet operation time and is given by the total time of all routes divided by the maximum repositioning time of the entire fleet.

\section {Algorithmic Proposal} 			\label{sec:Solution_Methodology}
To address our described SSBRP detailed in Section \ref{sec:Problem_Description}, we introduce a matheuristic procedure that combines a randomized multi-start algorithm with an integrated integer programming model. Our matheuristic, which relies on the \emph{Randomized Multi-Start Algorithm} (RMS), comprises two phases presented in Algorithm \ref{alg.Matheuristic}). 

In the first phase, a solution is constructed using RMS in the first phase, where the RMS creates a solution at each iteration, which is updated if it improves. The routes are created for each vehicle sequentially, i.e., one after the other. The routes are formed by iteratively inserting a new node at the end of the partial route using a greedy strategy, as detailed in Section \ref{sec:Randomized_Multi-Start_Algorithm}. 

The second phase involves optimizing the loading policy (or loading instructions) for a given feasible solution obtained by RMS in each iteration. This optimization is done through an integer mathematical model described in Section \ref{sec:Optimal_Loading_Instructions}. The objective is to determine the optimal loading policy for a given set of routes, specifying the number of operative and damaged bikes to be moved at the station or depot. The procedure is repeated in each iteration, yielding high-quality solutions within short computational times.
\begin{algorithm}[!htbp]
	\caption{-- \texttt{Matheuristic two-phase algorithm}} 		\label{alg.Matheuristic}\footnotesize 
	\begin{algorithmic}[1]
		\REQUIRE Instance $I$
		\REQUIRE Parameter $MaxIter$
		\STATE $\; Iter \gets 1$, $\; \mathcal{S} \gets \emptyset$
		\REPEAT
		\STATE $\quad \mathcal{S}^{'} \gets$Randomized Multi-Start$(I)\qquad\qquad\quad$\COMMENT{\textbf{Phase I}}	
		\STATE $\quad \mathcal{S}^{''} \gets$Optimal Loading Instructions$(\mathcal{S}^{'})\qquad$\COMMENT{\textbf{Phase II}}	
		\IF{$\mathcal{S}^{''}$ \textit{is better than} $\mathcal{S}$}
		\STATE $\mathcal{S} \gets \mathcal{S}^{''}$
		\STATE $Iter \gets 1$
		\ELSE					
		\STATE $Iter \gets Iter+1$
		\ENDIF										
		\UNTIL{$Iter=MaxIter$}				
		\ENSURE a feasible solution $\mathcal{S}$
	\end{algorithmic}
\end{algorithm}

In the subsequent sections, each phase of the developed matheuristic is presented in detail.

\section {Phase I: Randomized Multi-Start Algorithm}	\label{sec:Randomized_Multi-Start_Algorithm}
The \emph{Randomized Multi-start Algorithm} (RMS) sequentially constructs routes for each vehicle. Routes are formed by iteratively inserting a new node at the end of the partial route using a greedy strategy. Starting from the last visit $u$ of a partial route (or initially, the depot), we identify the set $F \subseteq V_0$ of feasible successors. For each candidate $v \in F$, we compute a ratio between the maximum number of bikes to be moved at $v$ and the travel time $t_{uv}$ we will take to visit each node candidate. Among the candidates with the highest ratio values, one node is randomly selected for insertion as the next visit in the current route. See Algorithm \ref{alg.RMS}.

At each algorithm stage, we have built a partial solution that includes station visits and loading and unloading operations on them. Suppose that the constructed route is $l$, and $u$ is the last node inserted in this partial route. As a result of the partial solution so far constructed, the states of the depot and the stations have changed (i.e., after each insertion of a new node $u$ at the end of the current route, say $l$, we update the whole information of the depot, stations, and vehicles). 

For each station $v$, let us denote by $\bar{d}_v$ and $\bar{a}_v$, respectively, the imbalance and the number of damaged bikes according to the partial solution built, and $\bar{p}^l_0$ the remaining number of operative bikes at the depot that can be taken by vehicle $l$ (that is, they have not been left by another vehicle). We denote by $\bar{V}_{def} = \{ \, v \in V \, : \, \bar{d}_v < 0 \, \}$, $\bar{V}_{spl} = \{ \, v \in V \, : \, \bar{d}_v > 0 \, \}$, $\bar{V}_{bal} = \{ \, v \in V \, : \, \bar{d}_v = 0 \, \}$ and $\bar{V}_{dam} = \{ \, v \in V \, : \, \bar{a}_v > 0 \, \}$. We denoted by $\bar{p}^l$ and $\bar{a}^l$, respectively, the number of operative bikes and damaged bikes on the vehicle $l$ after visiting $u$, and let $\bar{t}^l$ the partial traveling time of route $l$. Additionally, we denote by $\bar{k}^l$ the minimum number of free lockers in the vehicle in the section of the current route from the last visit to the depot to node $u$.

\begin{algorithm}[!htbp]
	\caption{-- \texttt{Randomized Multi-Start Algorithm}} 		\label{alg.RMS}\footnotesize 
	\begin{algorithmic}[1]
		\REQUIRE Instance $I$
		\STATE $\; \mathcal{S} \gets \emptyset$
		\FORALL{$l\in L$}
		\STATE $\; i \gets 1$, $\; F \gets \emptyset \qquad$\COMMENT{First stop}
		\STATE $ $ Initialize the route with $r^l(1)$ and Loading Instructions $y^l(1)$
		\STATE $ $ Define the set of candidate nodes $F \subseteq V_0$ to visit from $i$
		\WHILE{$|F|\neq \emptyset$}
		\STATE $ $ Calculate the ratios $\rho_v$ of each candidate $v \in F$
		\STATE $ $ Select randomly $v^*$ among the highest ratios: $\epsilon\times \rho_{max}  \leq \rho_{v} \leq \rho_{max}$
		\STATE $\; i \gets i+1 \qquad$\COMMENT{Next stop}
		\STATE $ $ Insert the node $v^*$ in the visit $r^l(i)$
		\STATE $ $ Define Loading Instructions $y^l(i)$ on the node $v^*$
		\STATE $ $ Define the set of new candidate nodes $F \subseteq V_0$ to visit from $i$
		\ENDWHILE
		\STATE $ $ Close the route $r^l$ visiting the depot	and delivery all the bikes		
		\STATE $ $ Update the route set $r^l$ and the Loading Instructions $y^l$
		\STATE $\; \mathcal{S} \gets (r^l, y^l)$		
		\ENDFOR		
		\ENSURE a feasible solution $\mathcal{S}$
	\end{algorithmic}
\end{algorithm}

The election of the next route visit $l$ is as follows: Firstly, a set of potential candidates $F$ is built. This set initialized with all the nodes $v \in \bar{V}_{def} \cup \bar{V}_{spl} \cup \bar{V}_{dam}$, $v \neq u$, such that $\bar{t}^l + t_{uv} + t_{v0} \leq T$, that is, those nodes not balanced or with damaged bikes that could be inserted after $u$ in the current route without exceeding the time limit $T$ (also considering the time needed to return to the depot). The depot is added to $F$ if $u \neq 0$ and $\bar{a}^l > 0$.

For each potential candidate $v \in F$, $v \neq 0$, we compute the maximum number of operative bikes $\beta_v$ and damaged bikes $\alpha_v$ we can move if $v$ is inserted after node $u$ in the current route. These quantities are computed as:

\begin{equation}  \label{eq:b_candidate}
	\beta_v = \left\{ \begin{array}{lcl} \mbox{min}\{ \, k^l - \bar{p}^l - \bar{a}^l, \, \bar{d}_v \, \} & & v \in F \cap \left( \bar{V}_{spl} \cup \bar{V}_{bal} \right) \\[0.2cm]
		
		\mbox{min}\{ \, \bar{p}^l+\mbox{min}\{\bar{p}^l_0, \bar{k}^l\}, \, |\bar{d}_v| \, \} & & v \in F \cap \bar{V}_{def} \\ \end{array} \right.
\end{equation}
\begin{equation}  \label{eq:a_candidate}
	\alpha_v = \left\{ \begin{array}{lcl} \mbox{min}\{ \, k^l - \bar{p}^l - \bar{a}^l - \beta_v, \, \bar{a}_v \, \} & & v \in F \cap \left( \bar{V}_{spl} \cup \bar{V}_{bal} \right) \\[0.2cm]
		
		\mbox{min}\{ \, k^l - \bar{p}^l - \bar{a}^l +\beta_v, \, \bar{a}_v \, \} & & v \in F \cap \bar{V}_{def} \\ \end{array} \right.
\end{equation}

Note that in the case where $v \in \bar{V}_{def}$ with $|\bar{d}_v| > \bar{p}^l$, we could consider to have taken more operative bikes from the depot if $\bar{p}^l_0 > 0$. To do it, we also have to check if the vehicle has enough capacity to carry it from the depot to $v$, represented by $\bar{k}^l$. This is the reason to consider $\mbox{min}\{ \, \bar{p}^l+\mbox{min}\{\bar{p}^l_0, \bar{k}^l\}, \, |\bar{d}_v| \, \}$ in (\ref{eq:b_candidate}). Finally, if $v=0$, we assume that the damaged bikes on the vehicle $\bar{a}^l$ would be unloaded in the depot.

The ultimate set of candidate successors for node $u$ in the route is derived by excluding from $F$ those stations $v$ where $\beta_v + \alpha_v = 0$. If $F = \emptyset$, we conclude the route by appending a visit to the depot and unloading all the bikes from the vehicle. Otherwise, for each $v \in F$, we calculate a ratio $\rho_v$ using the following formula:

\begin{equation}  \label{eq:ratio}
	\rho_v = \left\{ \begin{array}{lcl} \frac{(\alpha_v+\beta_v)^\theta}{t_{uv}}w_v & & v \in F \cap V \\[0.2cm]
		
		\frac{\bar{a}^l}{t_{u0}}\mu & & v = 0 \in F  \\ \end{array} \right.
\end{equation}

\noindent where, parameters $\theta \in (0, 1]$ and $\mu$ are employed to fine-tune the algorithm. A value of $\theta < 1$ is utilized to discourage the selection of stations with a substantial number of bikes to be moved, especially those that are far from $u$, in favor of stations that are closer to $u$ even if the number of bikes to be moved is not as large. Lastly, if $\mu > 1$, it promotes visits to the depot.

In a deterministic scheme, we would choose as the next visit in the route the candidate $v^*$ such that $\rho_{v^*} = \rho_{max}= \mbox{max}\{ \, \rho_v \, : \, v \in F \, \}$. However, we have embedded the algorithm into an randomized multi-start scheme that executes the construction of the solution a number of times equal to $MaxIter$ and output the best solution generated. 

In this algorithm, the next node $v^*$ is randomly selected among the candidates with the highest values of their ratios. Specifically, let $\epsilon$ be a random number in $(0,1)$, then $v^*$ is randomly selected among the candidates $v \in F$ such that $\epsilon\times \rho_{max}  \leq \rho_{v} \leq \rho_{max}$. 

Once the node $v^*$ has been selected for the next visit, we carry out the loading instructions procedure, that is, define the number of operative bikes and damaged bikes to be moved at the station or the depot are given by $\beta_{v^*}$ and $\alpha_{v^*}$, respectively. In the case where $|\bar{d}_{v^*}| > \bar{p}^l$ and $\mbox{min}\{\bar{p}^l_0, \bar{k}^l\} > 0$, $\beta_{v^*} - \bar{p}^l$ additional operative bikes are taken in the last visit to the depot, and the values of $\bar{k}^l$ and $\bar{p}^l_0$ are consequently updated.

\section {Phase II: Optimal Loading Instructions} 	\label{sec:Optimal_Loading_Instructions}
Given a feasible solution provided by RMS (See Section \ref{sec:Randomized_Multi-Start_Algorithm}) in each iteration described above, we apply an integer mathematical model to finding the optimal loading instructions in each route generated. Given a set of $s$ routes, this model determines the optimal loading policy, minimizing the final imbalance and the number of damaged bikes not removed. 

Below, we present the data and variables required for developing the proposed mathematical formulation.

\noindent \textbf{Data:}
\begin{itemize}
	\item initials of the stations and depot $v_\in V_o$ visited in $s$.
	\item from the fleet of vehicles.
		\begin{itemize}
			\item $\overline{L}$: Set of vehicles used in $s,\overline{L}\subseteq L$.
			\item $k^{l}$: Vehicle capacity $l,\overline{L}$.
		\end{itemize}		
	\item of the route.
		\begin{itemize}
			\item $n^{l}$: Number of visits of the route (vehicle) $l^{l}$.
			\item $r^{l}(i)$: Station or depot at the i-th visit of the route $l^{l}(i)$: Station or depot at the i-th visit of the route $l^{l}, i=1,\ldots,n^{l}$.
		\end{itemize}		
\end{itemize}

\noindent \textbf{Variables:}

\noindent  For each path $l\in \overline{L}$:
\begin{itemize}
	\item $b_o^{l}$: number of operative bikes in the depot that will use $l$.
	\item For each visit $i,i=1,\ldots,n^{l}$, to the stations or the depot:
	\begin{itemize}
		\item $x_i^{l}$: number of operative bikes taken (+) or left (-).
		\item $y_i^{l}$: number of damaged bikes taken (+) or left (-).
	\end{itemize}		
	\item For each visit $i,i=1,\ldots,n^{l}-1$, in vehicle $l$:
	\begin{itemize}
		\item $z_i^{l}$: number of operative bikes in $l$ when leaving the i-th visit.
		\item $w_i^{l}$: number of damaged bikes in $l$ when leaving the i-th visit.
	\end{itemize}				
\end{itemize}

\noindent \textbf{Objective Function:}
\begin{equation}\label{eq:fo}
	\text{Min }
	\sum_{v\in V_{spl}}(d_v-\sum_{l\in \overline{L}}\sum_{i \in I^l_v}x_i^{l})-
	\sum_{v\in V_{def}}(d_v-\sum_{l\in \overline{L}}\sum_{i \in I^l_v}x_i^{l})+
	\sum_{v\in V_{dam}}(a_v-\sum_{l\in \overline{L}}\sum_{i \in I^l_v}y_i^{l})	
\end{equation}
\textbf{s.t.:}\
\begin{equation}\label{eq:cap_veh}
	\sum_{i=1}^{j}(x_{i}^{l}+y_{i}^{l})\leq k^{l} 	
	\qquad  \forall l\in \overline{L}, j=1,\ldots,n^{l}-1
\end{equation}
\begin{equation}\label{eq:pos_veh}
	\sum_{i=1}^{j}x_{i}^{l}\geq 0 	
	\qquad \forall l\in \overline{L}, j=1,\ldots,n^{l}-1
\end{equation}
\begin{equation}\label{eq:descarga_buenas_veh}
	x_{n^{l}}^{l} = -\sum_{i=1}^{n^{l}-1}x_{i}^{l}
	\qquad \forall l\in \overline{L}
\end{equation}
\begin{equation}\label{eq:descarga_ave_veh}
	\sum_{i=1}^{j}y_{i}^{l} = 0	
	\qquad \forall l\in \overline{L}, \forall j \in I^l_0
\end{equation}
\begin{equation}\label{eq:R4}
	\sum_{l\in \overline{L}} \sum_{i \in I^l_v}x^{l}_{i} \leq d_{v}
	\qquad \forall v\in V_{spl}
\end{equation}
\begin{equation}\label{eq:R5}
	\sum_{l\in \overline{L}} \sum_{{i \in I^l_v}}x^{l}_{i} \geq d_{v}
	\qquad \forall v\in V_{def}
\end{equation}
\begin{equation}\label{eq:R6}
	\sum_{l\in \overline{L}} \sum_{{i \in I^l_v}}y^{l}_{i} \leq a_{v}
	\qquad \forall v\in V_{dam}
\end{equation}
\begin{equation}\label{eq:R7}
	\sum_{l\in \overline{L}}w^{l}_{0}\leq p_o
\end{equation}
\begin{equation}\label{eq:depot}
	x^{l}_{j} \leq w^{l}_{0} - \sum_{{i \in I^l_0, i<j}} x^{l}_{i}
	\qquad \forall l\in \overline{L}, \forall {j \in I^l_0}	
\end{equation}
\begin{equation}\label{eq:R18a}
	{x^{l}_{i} = 0
		\qquad \forall l\in \overline{L}, \forall v\in V_{bal}, \forall i \in I^l_v}	
\end{equation}
\begin{equation}\label{eq:R18b}
	{y^{l}_{i} = 0
		\qquad \forall l\in \overline{L}, \forall v\in V \backslash V_{dam}, \forall i \in I^l_v}	
\end{equation}
\begin{equation}\label{eq:R18}
	x^{l}_{i} \in \mathbb{Z_+}
	\qquad \forall l\in \overline{L}, \forall {v\in V_{spl}, \forall i \in I^l_v}	
\end{equation}
\begin{equation}\label{eq:R19}
	x^{l}_{i} \in \mathbb{Z_+}
	\qquad \forall l\in \overline{L}, \forall {v\in V_{def}, \forall i \in I^l_v}
\end{equation}
\begin{equation}\label{eq:R20}
	y^{l}_{i} \in \mathbb{Z_+}
	\qquad \forall l\in \overline{L}, \forall {v\in V_{dam}, \forall i \in I^l_v}	
\end{equation}
\begin{equation}\label{eq:R20a}
	{y^{l}_{i} \in \mathbb{Z_+}
		\qquad \forall l\in \overline{L}, \forall i \in I^l_0}	
\end{equation}
\begin{equation}\label{eq:var_w}
	w^{l}_0 \in \mathbb{Z_+}
	\qquad \forall l\in \overline{L}
\end{equation}
\begin{equation}\label{eq:R22}
	x^{l}_{i} \in \mathbb{R}
	\qquad \forall l\in \overline{L}, {\forall i \in I^l_0}
\end{equation}

The objective function minimizes the total imbalance of the stations and the number of damaged bikes that have not been collected from the stations. Constraints (\ref{eq:cap_veh}-\ref{eq:pos_veh}) guarantee that the number of operative bikes and the total number of bikes on the vehicle is non-negative and below capacity, respectively, at each visit. Constraints (\ref{eq:descarga_buenas_veh}) ensure that all operative bikes are unloaded at the end of each route. Constraints (\ref{eq:descarga_ave_veh}) guarantee that all damaged bikes are unloaded at each visit to the depot. Constraints (\ref{eq:R4}-\ref{eq:R6}) limit the number of bikes that can be managed in the stations. Constraints (\ref{eq:R7}) limit the number of operative bikes that can be taken from the depot by all the vehicles, and constraints (\ref{eq:depot}) control the number of operative bikes that can be taken in each visit to the depot. Finally, constraints (\ref{eq:R18a}-\ref{eq:R22}) define the domain of the variables.


\section {Computational Experiments} 		\label{sec:Computational_Experiments}
The computational analysis to assess the performance of our algorithm is structured as follows: We describe the instances based on real-world data used for testing and introduce the parameters information used in our solution method. Then, we evaluate the impact of the route construction parameters on the best solutions found. Finally, we present a comparative study with another SBRP variant that has been more extensively studied in the literature.

We evaluated the performance of our algorithms using 136 instances, categorized into two groups: Palma and Wien.

The first set comprises instances from the Bike-Sharing System in Palma de Mallorca, Spain, as initially proposed in \cite{Alvarez-Valdes2016}. Each instance includes 28 stations and one depot, with an initial inventory of 10 operative bikes, represented as $p_0=10$. Additionally, two fleet sizes (2 and 3 vehicles with a capacity of $k^l=20$) and two variations of maximum time (2 and 4 hours) were considered. This set consists of two instances in each group, totaling 56 instances. 	

The second set was adapted from \cite{Rainer-Harbach2015}, with modifications to align with our problem by randomly substituting some operational bikes in the stations with damaged bikes. These instances involve 20, 30, 60, and 90 stations and one depot (with no initial inventory, i.e., $p_0=0$). For this set, two values for the maximum time (4 and 8 hours) and three fleet sizes (2, 3, and 5 vehicles with a capacity of $k^l=20$) were considered. Each group has five instances, making a total of 80 instances in this set.

All algorithms were implemented in C++ and executed on a PC with an Intel(R) Core(TM) i5-4200 CPU, operating at 1.60 GHz, and equipped with 4.00 GB of RAM. Throughout the experiments, only a single thread was utilized. The integer programming model described in Section \ref{sec:Optimal_Loading_Instructions} was resolved using IBM ILOG CPLEX 20.0. 

Our experimentation evaluated the parameters: $MaxIter=500$, $\theta$ with values 0.3, 0.5, and 0.8, and $\mu$ with values 1.0, 1.5, and 2.0. The scaling factors in the objective function described in Section \ref{sec:Problem_Description} were set to $\gamma_d = \gamma_a = \gamma_t = 1$. By employing these factors, enhancing the system balance and collecting damaged bikes will consistently have a greater impact on the objective value than reducing the total time of the routes.

In our experimentation, first we independently assessed the impact of the route construction parameters $\theta$ and $\mu$ on the best solutions found.

Table \ref{Table.ConstructionParam} displays the average results for each group of cases. The first two columns show the route construction parameters $\theta$ and $\mu$, respectively. The subsequent columns present the instance sets from Palma and Wien, showing the average values of the objective function ($O.F$), the average number of iterations needed to achieve the best solution ($Iter$), and the time of ($CPU$) in seconds required to obtain those solutions.

\begin{table}[H]
	\caption{Average results of the route construction parameters $\theta$ and $\mu$.} \label{Table.ConstructionParam}
	\begin{center}
		\scriptsize 
		\begin{tabular}{|c|c|ccc|ccc|}
			\hline			
			$\theta$ & $\mu$ & $O.F_{Palma}$ & $Iter_{Palma}$ & $CPU_{Palma}$ & $O.F_{Wien}$ & $Iter_{Wien}$ & $CPU_{Wien}$\\
			\hline			
			\textbf{0.3} & \textbf{1} & \textbf{0.745} & 133.21 & 1.65 & 1.010 & 229.56 & 1.14 \\
			\textbf{0.3} & \textbf{1.5} & \textbf{0.745} & 133.21 & 3.22 & 1.010 & 235.24 & 1.62 \\
			\textbf{0.3} & \textbf{2} & \textbf{0.745} & 133.21 & 3.30 & 1.011 & 232.61 & 0.97 \\
			\hline
			0.5 & 1 & 0.748 & 152.07 & 3.34 & 0.998 & 238.09 & 2.36 \\
			0.5 & 1.5 & 0.748 & 147.79 & 1.32 & 0.998 & 240.24 & 3.05 \\
			0.5 & 2 & 0.748 & 147.79 & 4.09 & 0.997 & 249.09 & 3.07 \\
			\hline
			0.8 & 1 & 0.755 & 129.80 & 2.53 & 0.992 & 212.25 & 4.33 \\
			0.8 & 1.5 & 0.755 & 129.80 & 3.57 & 0.992 & 213.28 & 3.63 \\
			\textbf{0.8} & \textbf{2} & 0.755 & 129.80 & 1.12 & \textbf{0.980} & 215.90 & 4.27 \\
			\hline
		\end{tabular}
	\end{center}
\end{table}

Table \ref{Table.ConstructionParam} displays the outcomes, revealing that, on average, the parameter $\theta=0.3$ consistently yields better solutions for the Palma instances, irrespective of the $\mu$ value used. On the other hand, for the Wien instances, the optimal parameters for achieving the best solutions, on average, are $\theta=0.8$ and $\mu=2$. Subsequently, using these identified best parameters, we evaluated the 136 instances. 

For the Palma instances, our method successfully achieved complete station balancing in 96.9\% of cases and efficiently collected all damaged bikes. The average CPU time for this set was 2.09 seconds. Similar favorable outcomes were observed for instances in Wien with 20 stations, where achieved complete station balancing in 95.7\%. However, for instances with 30-90 stations in Wien posed more challenges. In these cases, nearly all vehicles were utilized, and in many instances, the solution could only fulfill some tasks, specifically achieving complete station balancing and retrieving all damaged bikes. We achieved complete station balancing in 62.35\%. Notably, improvements were observed with an increase in the number of vehicles or the maximum time $T$. The average CPU time required for this group (30-90 stations) was 2.65 seconds.

Finally, we compare the performance of our solution methodology with a variant presented in the study by Rainer-Harbach et al.~\cite{Rainer-Harbach2015}. These instances are based on the Bike-Sharing System (BSS) of Vienna, Austria, and include 10 and 20 stations with one depot. For this set, three values for the maximum time (2, 4, and 8 hours) and three fleet sizes (1, 2, and 3 vehicles) were considered. Each instance set uses a unique combination of $|V|, |L|, T$ and contains 30 cases, resulting in a total of 12 sets and 360 instances.

To compare our approach with the variant proposed in \cite{Rainer-Harbach2015}, we need to align the characteristics of our problem with theirs. We believe that the features of our problem (described in Section \ref{sec:Problem_Description}) can be considered a general variant of the BRP compared to those found in the literature. Thus, our framework can accommodate the study of other BRP variants, including the one proposed in \cite{Rainer-Harbach2015}. To do this, we may need to simplify or relax specific characteristics or assumptions of our problem, particularly regarding the nodes, the routes, and the objective function.

After adapting our algorithm to the variant proposed in \cite{Rainer-Harbach2015}, we apply it to the original instances mentioned earlier. Table \ref{Table.ResulRainer} presents the average results for each group of cases. The table first displays the results of the MIP model from \cite{Rainer-Harbach2015}, including the mean upper bounds ($\overline{ub}$), the mean lower bounds ($\overline{lb}$) and the median total run times ($t_{MIP}(s)$). Second, it shows the best-known solutions ($\overline{bks}$) obtained on average by Rainer-Harbach et al.~\cite{Rainer-Harbach2015}. Third, it presents the average values of the objective function ($\overline{obj}$) generated by our proposed solution methodology from Section \ref{sec:Solution_Methodology}. Fourth, it indicates the gap between $\overline{obj}$ and the mean lower bounds $\overline{lb}$ obtained in the MIP model. Finally, it demonstrates the difference between $\overline{obj}$ and $\overline{bks}$ and the mean runtime ($CPU_{obj}(s)$).

\begin{table}[!htbp]
	\caption{Average results from 12 instances for each group of $|V|\in \{10, 20\}$.} \label{Table.ResulRainer}
	\begin{center}
		\scriptsize 
		\begin{tabular}{|c|c|c|c|c|c|c|c|c|c|c|}
			\hline
			$|V|$ & $|L|$ & $T(h)$ & $\overline{ub}$ & $\overline{lb}$ & $t_{MIP}(s)$ & $\overline{bks}$ & $\overline{obj}$ & gap & dif & $CPU_{obj}(s)$\\
			\hline
			10 & 1 & 2 & 27.80143  & 27.80143  & 3.0  & 27.86810  & 28.00143  & \textbf{0.00719}  & 0.13333  & 0.22\\
			10 & 1 & 4 & 3.46948  & 0.17536  & 3.600  & 3.50949  & 3.53615  & 19.16510  & \textbf{0.02666}  & 0.29\\
			10 & 1 & 8 & 0.00322  & 0.00260  & 3.600  & 0.00319  & 0.00320  & 0.22932  & \textbf{0.00000}  & 0.18\\
			10 & 2 & 2 & 9.13607  & 8.80436  & 897,5  & 9.26937  & 9.32270  & \textbf{0.05887}  & \textbf{0.05333}  & 1.38\\
			10 & 2 & 4 & 0.00341  & 0.00326  & 3.600  & 0.00339  & 0.00339  & \textbf{0.04047}  & \textbf{0.00001}  & 0.78\\
			10 & 2 & 8 & 0.00324  & 0.00315  & 3.600  & 0.00319  & 0.00320  & \textbf{0.01396}  & \textbf{0.00000}  & 0.70\\
			\hline
			20 & 2 & 2 & 52.46973  & 29.77982  & 3.600  & 51.52307  & 51.64306  & 0.73416  & 0.11999  & 0.46\\
			20 & 2 & 4 & 16.00568  & 0.00387  & 3.600  & 5.24584  & 5.29917  & 1367.23341  & \textbf{0.05333}  & 0.28\\
			20 & 2 & 8 & 0.14154  & 0.00350  & 3.600  & 0.00638  & 0.00638  & 0.82338  & \textbf{0.00000}  & 0.59\\
			20 & 3 & 2 & 34.73765  & 2.08780  & 3.600  & 28.75112  & 28.96445  & 12.87320  & 0.21333  & 1.73\\
			20 & 3 & 4 & 0.94106  & 0.00536  & 3.600  & 0.00670  & 0.00670  & 0.25089  & \textbf{0.00001}  & 1.54\\
			20 & 3 & 8 & 8.14200  & 0.00344  & 3.600  & 0.00638  & 0.00638  & 0.85440  & \textbf{0.00001} & 1.13\\					
			\hline
		\end{tabular}
	\end{center}
\end{table}

For the instance groups where the MIP showed small gaps between upper and lower bounds, our solution method found solutions with equal or slightly different objective values. 

The results in the $\overline{bks}$ column of Table \ref{Table.ResulRainer} represent the solutions, on average, with the best objective values among the four variants proposed in \cite{Rainer-Harbach2015}. We observe a clear tendency that our solution method performs similarly to the variants proposed in \cite{Rainer-Harbach2015} in most instances.

It is important to note that the variants proposed by \cite{Rainer-Harbach2015} are specifically designed to solve their problem. At the same time, our approach needs to be adapted beforehand because our solution methods were initially developed for solving the BRP variant described in Section \ref{sec:Problem_Description}. Nevertheless, we achieved similar results in 83.3\% of the instances tested.

\section {Concluding Remarks}  				\label{sec:Concluding_Remarks}
This paper introduces a specific instance of the \textit{Static Bike-sharing Repositioning Problem} (SBRP). It incorporates new considerations for existing SBRP to improve their realism, such as operational and damaged bikes, a heterogeneous fleet, and multiple visits between stations and the depot. 

To solve this problem, we propose a matheuristic approach based on a randomized multi-start algorithm integrated with integer programming, incorporating various strategies and route construction parameters that, in combination, contribute to improving the quality of the solutions found.

To our knowledge, there are currently no previous works that allow us to compare the quality of the solutions obtained. However, we are developing a method to calculate a lower bound that allows us to evaluate the quality of the solutions.

To demonstrate the applicability of our approach, we adapted our algorithm to a more straightforward problem studied in the literature, achieving competitive results. This suggests that the characteristics of our problem are a general variant of the BRP compared to those described in existing literature.

We proposed an integer programming model for optimizing the number of operative and damaged bikes that will be moved between stations and/or the depot (loading instructions). This model can be solved very quickly with IBM ILOG CPLEX for the size of the instances we have tried.

\begin{credits}

\end{credits}

%
%
%
\bibliographystyle{splncs04}
\bibliography{mybibliography2024}
%

%
%
%
%
\end{document}